\documentclass[12pt]{article}
\usepackage{amssymb}
\usepackage{amsfonts}
\usepackage{amsmath}
\usepackage{amsbsy}
\newcommand{\be}{\begin{equation}}
\newcommand{\ee}{\end{equation}}
\newcommand{\bea}{\begin{eqnarray}}
\newcommand{\eea}{\end{eqnarray}}
\newcommand{\ba}{\begin{array}}
\newcommand{\ea}{\end{array}}
\newcommand{\R}{I\!\!R}

\newcommand{\tp}{I\!\!L}

\newcommand{\bc}{\begin{center}}
\newcommand{\ec}{\end{center}}
\newcommand{\ben}{\begin{enumerate}}
\newcommand{\een}{\end{enumerate}}
\newcommand{\bfi}{\begin{figure}}
\newcommand{\efi}{\end{figure}}

\newcommand{\bq}{\begin{quote}}
\newcommand{\eq}{\end{quote}}
\newcommand{\bqu}{\begin{quotation}}
\newcommand{\equ}{\end{quotation}}
\newenvironment{emphit}{\begin{itemize}}{\end{itemize}}
\newcommand{\bemp}{\begin{emphit}}
\newcommand{\eemp}{\end{emphit}}

\newcommand{\bt}{\begin{tabular}}
\newcommand{\et}{\end{tabular}}

\newtheorem{myth}{Theorem}[section]
\newtheorem{mylem}{Lemma}[section]
\newtheorem{mycor}{Corollary}[section]

\newtheorem{mydef}{Definition}[section]
\newtheorem{myrem}{Remark}[section]

\begin{document}
\date{}
\title{Some necessary and sufficient conditions for parastrophic invariance of the associative law
in quasigroups
}
\author{T\`em\'it\'op\'e Gb\'ol\'ah\`an Ja\'iy\'e\d ol\'a\\Department of
Mathematics,\\
Obafemi Awolowo University, Ile Ife,
Nigeria.\\jaiyeolatemitope@yahoo.com,~tjayeola@oauife.edu.ng}
\maketitle

\begin{abstract}
Every quasigroup $(S,\cdot )$ belongs to a set of 6 quasigroups,
called parastrophes denoted by $(S,\pi_i)$, $i\in \{1,2,3,4,5,6\}$.
It is shown that isotopy-isomorphy is a necessary and sufficient
condition for any two distinct quasigroups $(S,\pi_i)$ and
$(S,\pi_j)$, $i,j\in \{1,2,3,4,5,6\}$ to be parastrophic invariant
relative to the associative law. In addition, a necessary and
sufficient condition for any two distinct quasigroups $(S,\pi_i)$
and $(S,\pi_j)$, $i,j\in \{1,2,3,4,5,6\}$ to be parastrophic
invariant under the associative law is either if the
$\pi_i$-parastrophe of $H$ is equivalent to the $\pi_i$-parastrophe
of the holomorph of the $\pi_i$-parastrophe of $S$ or if the
$\pi_i$-parastrophe of $H$ is equivalent to the $\pi_k$-parastrophe
of the $\pi_i$-parastrophe of the holomorph of the
$\pi_i$-parastrophe of $S$, for a particular $k\in \{1,2,3,4,5,6\}$.
\end{abstract}

\section{Introduction}
\paragraph{}
Let $L$ be a non-empty set. Define a binary operation ($\cdot $) on
$L$ such that $x\cdot y\in L$ for all $x,y$ in $L$. Then, $(L, \cdot
)$ is called a groupoid. If the system of equations ; $a\cdot x=b$
and $y\cdot a=b$ have unique solutions for $x$ and $y$ respectively,
then $(L, \cdot )$ is called a quasigroup. At times we write
$y=x\backslash z$ or equivalently $x=z/y$ for $x\cdot y=z$.
Furthermore, if there exists a unique element $e\in L$ called the
identity element such that for all $x$ in $L$, $x\cdot e=e\cdot
x=x$, $(L, \cdot )$ is called a loop. On the other hand if in the
quasigroup $(L,\cdot )$, $xy\cdot z=x\cdot yz$ for all $x,y,z$ in
$L$(associativity property), then $(L,\cdot )$ is called a group and
has an identity element(\cite{phd2}). So, a loop is not a group but
a group is a loop. This is the foundational reason for the study of
non-associative algebraic system

It has been noted that every quasigroup $(L,\cdot )$ belongs to a
set of 6 quasigroups, called adjugates by Fisher and Yates
\cite{phd72}, conjugates by Stein \cite{phd73}, \cite{phd76} and
Belousov \cite{phd77} and parastrophes by Sade \cite{phd74}. They
have been studied by Artzy \cite{phd71}, Charles Lindner and Dwight
Steedley \cite{phd173.1} and a detailed study on them can be found
in \cite{phd3}, \cite{phd39} and \cite{phd49}. The most recent study
of the parastrophes of a quasigroup(loop) are by Sokhatskii
\cite{phd171,phd172}, Duplak \cite{phd154} and Shchukin and Gushan
\cite{phd173}. For a quasigroup $(L,\cdot )$, its parastrophes are
denoted by $(L,\pi_i)$, $i\in \{1,2,3,4,5,6\}$ hence one can take
$(L,\cdot )=(L,\pi_1)$. A quasigroup which is equivalent to all its
parastrophes is called a totally symmetric quasigroup(introduced by
Bruck \cite{phd165}) while its loop is called a Steiner loop. For
more on quasigroup, loops and their properties, readers should check
\cite{phd3}, \cite{phd41}, \cite{phd39}, \cite{phd49}, \cite{phd42}
and \cite{phd75}. Let $(G,\oplus)$ and $(H,\otimes)$ be two distinct
quasigroups. The triple $(A,B,C)$ such that
\begin{displaymath}
A,B,C~:~(G,\oplus)\longrightarrow (H,\otimes)
\end{displaymath}
are bijections is said to be an isotopism if and only if
\begin{displaymath}
xA\otimes yB=(x\oplus y)C~\forall~x,y\in G.
\end{displaymath}
Thus, $H$ is called an isotope of $G$ and they are said to be
isotopic. According to \cite{phd3}, every group is a G-loop(i.e a
loop that is isomorphic to all its loop isotopes). Hence, every loop
isotope of a group is a group but this is not true for quasigroup
isotopes of a group, they are not necessarily associative. The aim
of this work is to find necessary and sufficient condition(s) for
any parastrophe of a group to be a group.

Let $Aum(L,\theta)$ be the automorphism group of a loop(quasigroup)
$(L,\theta)$, and the set $H=Aum(L,\theta)\times (L,\theta)$. If we
define '$\circ$' on $H$ such that
\begin{displaymath}
(\alpha, x)\circ (\beta, y)=(\alpha\beta, x\beta\theta
y)~\forall~(\alpha, x),(\beta, y)\in H,
\end{displaymath}
then $H(L,\theta)=(H,\circ)$ is a loop(quasigroup) according to
Bruck \cite{phd82} and is called the Holomorph of $(L,\theta )$. It
can be observed that a loop is a group if and only if its holomorph
is a group.

Interestingly, Adeniran \cite{phd79} and Robinson \cite{phd85},
Oyebo and Adeniran \cite{phdoyebo}, Chiboka and Solarin
\cite{phd80}, Bruck \cite{phd82}, Bruck and Paige \cite{phd40},
Robinson \cite{phd7} and Huthnance \cite{phd44} have respectively
studied the holomorphs of Bruck loops and Bol loops, central loops,
conjugacy closed loops, inverse property loops, A-loops, extra loops
and weak inverse
property loops.\\

In this paper, it is proved that $(S,\cdot )$ is an associative
quasigroup if and only if any one of four particular parastrophes of
$(S,\cdot )$ obeys a Khalil condition or Evans' generalized
associative law or Belousov's balanced identity or Falconer's
generalized group identity. It is shown that isotopy-isomorphy is a
necessary and sufficient condition for any two distinct quasigroups
$(S,\pi_i)$ and $(S,\pi_j)$, $i,j\in \{1,2,3,4,5,6\}$ to be
parastrophic invariant relative to the associative law.

Furthermore, the relationship between the parastrophes of the
holomorph of a quasigroup and the holomorphs of the parastrophes of
the same quasigroup are investigated. The following results are
proved. For a quasigroup $L$ with holomorph $H$, it is shown that
the $\pi_i$-parastrophe of $H$ is equivalent to the
$\pi_i$-parastrophe of the holomorph of the $\pi_i$-parastrophe of
$L$ if and only if $L$ is equivalent to its $\pi_i$-parastrophe for
each $i=1,2,3,4,5,6$. A necessary and sufficient condition for any
two distinct quasigroups $(S,\pi_i)$ and $(S,\pi_j)$, $i,j\in
\{1,2,3,4,5,6\}$ to be parastrophic invariant under the associative
law is either if the $\pi_i$-parastrophe of $H$ is equivalent to the
$\pi_i$-parastrophe of the holomorph of the $\pi_i$-parastrophe of
$S$ or if the $\pi_i$-parastrophe of $H$ is equivalent to the
$\pi_k$-parastrophe of the $\pi_i$-parastrophe of the holomorph of
the $\pi_i$-parastrophe of $S$, for a particular $k\in
\{1,2,3,4,5,6\}$.

\section{Preliminaries}
\begin{mydef}\label{par:loop}
Let $(L,\theta )$ be a quasigroup. The 5 parastrophes or conjugates
or adjugates of $(L,\theta )$ are quasigroups whose binary
operations
$\theta^*~,~\theta^{-1}~,~{}^{-1}\theta~,~(\theta^{-1})^*~,~({}^{-1}\theta
)^*$ defined on $L$ are given by :
\begin{description}
\item[(a)] $(L,\theta^*)~: ~y\theta ^*x=z\Leftrightarrow x\theta y=z~\forall~x,y,z\in
L$.
\item[(b)] $(L,\theta^{-1})~:~x\theta ^{-1}z=y\Leftrightarrow x\theta y=z~\forall~x,y,z\in
L$.
\item[(c)] $(L,{}^{-1}\theta )~:~z~{}^{-1}\theta y=x\Leftrightarrow x\theta y=z~\forall~x,y,z\in
L.$.
\item[(d)] $(L,(\theta ^{-1})^*)~:~z(\theta ^{-1})^*x=y\Leftrightarrow x\theta y=z~\forall~x,y,z\in
L$.
\item[(e)] $(L,({}^{-1}\theta )^*)~:~y({}^{-1}\theta )^*z=x\Leftrightarrow x\theta y=z~\forall~x,y,z\in
L$.
\end{description}
The five parastrophes of $H(L,\theta)=(H,\circ)$ shall be denoted by
\begin{displaymath}
(H,\circ^*),~(H,\circ^{-1}),~(H,{}^{-1}\circ),~\big(H,(\circ^{-1})^*,\big)~
\textrm{and}~\big(H,({}^{-1}\circ)^*\big).
\end{displaymath}

The holomorphes of the five parastrophes of the quasigroup
$(L,\theta)$ are denoted by :
\begin{displaymath}
H(L,\theta^*)=(H,\circ_*),~H(L,\theta^{-1})=(H,\circ_{-1}),~H(L,{}^{-1}\theta)=(H,{}_{-1}\circ),
\end{displaymath}
\begin{displaymath}
~H\big(L,(\theta^{-1})^*\big)=\big(H,(\circ_{-1})_*\big)~\textrm{and}~
H\big(L,({}^{-1}\theta)^*\big)=\big(H,{}_*({}_{-1}\circ)\big).
\end{displaymath}
\end{mydef}

\begin{mydef}\label{rep:loop}
Let $(L,\theta)$ be a quasigroup.
\begin{description}
\item[(a)] $R_x$ and $L_x$ represent the right and left translation maps in
$(L,\theta )~\forall~x\in L$.
\item[(b)] $R_x^*$ and $L_x^*$ represent the right and left translation maps in
$(L,\theta^*)~\forall~x\in L$.
\item[(c)] ${\cal R}_x$ and ${\cal L}_x$ represent the right and left translation maps in
$(L,\theta^{-1})~\forall~x\in L$.
\item[(d)] ${\R}_x$ and ${\tp}_x$ represent the right and left translation maps in
$(L,{}^{-1}\theta)~\forall~x\in L$.
\item[(e)] ${\cal R}_x^*$ and ${\cal L}_x^*$ represent the right and left translation maps in
$(L,(\theta^{-1})^*)~\forall~x\in L$.
\item[(f)] ${\R}_x^*$ and ${\tp}_x^*$ represent the right and left translation maps in
$(L,({}^{-1}\theta)^*)~\forall~x\in L$.
\end{description}
\end{mydef}

\begin{myrem}
If $(L,\theta )$ is a loop, $(L,\theta^*)$ is also a loop(and vice
versa) while the other adjugates are quasigroups. Furthermore,
$(L,\theta^{-1})$ and $(L,({}^{-1}\theta )^*)$ have left identity
elements, that is they are left loops while $(L,{}^{-1}\theta )$ and
$(L,(\theta ^{-1})^*)$ have right identity elements, that is they
are right loops. $(L,\theta^{-1})$ or $(L,{}^{-1}\theta )$ or
$(L,(\theta ^{-1})^*)$ or $(L,({}^{-1}\theta )^*)$ is a loop if and
only if $(L,\theta )$ is a loop of exponent $2$.
\end{myrem}

\begin{mylem}\label{0:1}
If $(L,\theta )$ is a quasigroup, then
\begin{enumerate}
\item $R_x^*=L_x~,~L_x^*=R_x~,~{\cal L}_x=L_x^{-1}~,~{\R}_x=R_x^{-1}~,~{\cal R}_x^*=L_x^{-1}~,~{\tp}_x^*=R_x^{-1}~\forall~x\in L$.
\item ${\cal L}_x=R_x^{*-1}~,~{\R}_x=L_x^{*-1}~,~{\cal R}_x^*=R_x^{*-1}={\cal L}_x~,~{\tp}_x^*=L_x^{*-1}={\R}_x~\forall~x\in L$.
\end{enumerate}
\end{mylem}
{\bf Proof}\\
The proof of these follows by using Definition~\ref{par:loop} and
Definition~\ref{rep:loop}.
\begin{enumerate}
\item $y\theta^*x=z\Leftrightarrow x\theta y=z\Rightarrow y\theta^*x=x\theta y
\Rightarrow yR_x^*=yL_x\Rightarrow R_x^*=L_x$. Also,
$y\theta^*x=x\theta y\Rightarrow
xL_y^*=xR_y\Rightarrow L_y^*=R_y$.\\
$x\theta^{-1}z=y\Leftrightarrow x\theta y=z\Rightarrow x\theta
(x\theta^{-1}z)=z \Rightarrow x\theta z{\cal L}_x=z\Rightarrow
z{\cal L}_xL_x=z\Rightarrow {\cal L}_xL_x=I$. Also,
$x\theta^{-1}(x\theta y)=y\Rightarrow x\theta^{-1}yL_x=y\Rightarrow
yL_x{\cal L}_x=y \Rightarrow L_x{\cal L}_x=I$.
Hence, ${\cal L}_x=L_x^{-1}~\forall~x\in L$.\\
$z({}^{-1}\theta ) y=x\Leftrightarrow x\theta y=z\Rightarrow
(x\theta y)({}^{-1}\theta ) y=x\Rightarrow xR_y({}^{-1}\theta )
y=x\Rightarrow xR_y{\R}_y=x\Rightarrow R_y{\R}_y=I$. Also,
$(z({}^{-1}\theta ) y)\theta y=z\Rightarrow z{\R}_y\theta
y=z\Rightarrow z{\R}_yR_y=z\Rightarrow
{\R}_yR_y=I$. Thence, ${\R}_y=R_y^{-1}~\forall~x\in L$.\\
$z(\theta^{-1})^*x=y\Leftrightarrow x\theta y=z$, so, $x\theta
(z(\theta^{-1})^*x)=z\Rightarrow x\theta z{\cal R}_x^*=z \Rightarrow
z{\cal R}_x^*L_x=z\Rightarrow {\cal R}_x^*L_x=I$. Also, $(x\theta
y)(\theta^{-1})^*x=y\Rightarrow yL_x(\theta ^{-1})^*x=y\Rightarrow
yL_x{\cal R}_x^*=y\Rightarrow L_y{\cal R}_x^*=I$.
Whence, ${\cal R}_x^*=L_x^{-1}$.\\
$y({}^{-1}\theta )^*z=x\Leftrightarrow x\theta y=z$, so,
$y({}^{-1}\theta )^*(x\theta y)=x\Rightarrow
y({}^{-1}\theta)^*xR_y=x\Rightarrow xR_y{\tp}_y^*=x\Rightarrow
R_y{\tp}_y^*=I$. Also, $(y({}^{-1}\theta)^*z)\theta y=z\Rightarrow
z{\tp}_y^*\theta y=z\Rightarrow z{\tp}_y^*R_y=z\Rightarrow
{\tp}_y^*R_y=I$. Thus, ${\tp}_y^*=R_y^{-1}$.
\item These ones follow from 1..
\end{enumerate}

\begin{myth}\label{0:1.3}(Falconer Theorem~2.9 \cite{phd159})

If a quasigroup $Q$ is isotopic to a group $G$, then all the
parastrophes of $Q$ are isotopic to $G$.
\end{myth}

In the past, some isotopy closure properties for groups(i.e
necessary and sufficient conditions for a quasigroup to be isotopic
to a group) and isotopy-isomorphy conditions for groups(i.e
necessary and sufficient conditions for the isomorphism of
quasigroups isotopic to a group \cite{phd181}) have been proved.
Some of the formal are stated below.

\begin{myth}\label{0:1.4}(Evans \cite{phd167})

A quasigroup $(Q,\cdot )$ is isotopic to a group if and only if $Q$
obeys the law
\begin{displaymath}
[(xP_1\cdot yP_2)P_3\cdot zP_4]P_5=[xQ_1\cdot (yQ_2\cdot
zQ_3)Q_4]Q_5\qquad\textrm{Evans law}
\end{displaymath}
where $P_i~,Q_i,~i=1,2,3,4,5$ are permutations on $(Q,\cdot )$.
\end{myth}

\begin{myth}\label{0:1.5}(Belousov \cite{phd166})

A quasigroup $(Q,\cdot )$ is isotopic to a group if and only if $Q$
obeys an identity
\begin{displaymath}
w_1(x_1,x_2,\cdots x_n)=w_2(x_1,x_2,\cdots
x_n)\qquad\textrm{balanced identity}
\end{displaymath}
where each $x_i$, $i=1,2,\cdots ,n$, occurs exactly once in $w_1$
and in $w_2$ and $w_1$ and $w_2$ involve different grouping of some
triple.
\end{myth}

\begin{myth}\label{0:1.6}(Falconer Theorem~2.10 \cite{phd159})

A quasigroup $(Q,\cdot )$ is isotopic to a group if and only if $Q$
satisfies a generalized group identity
\end{myth}

\begin{myth}\label{0:1.2}(Khalil Conditions \cite{phd87})

A quasigroup is an isotope of a group if and only if any one of the
following six identities are true in the quasigroup for all elements
$x,y,z,u,v$.
\begin{enumerate}
\item $x\{z\backslash [(z/u)v]\}=\{[x(z\backslash z)]/u\}v$
\item $x\{u\backslash [(z/u)v]\}=\{[x(u\backslash z)]/u\}v$
\item $x\{z\backslash [(u/u)v]\}=\{[x(z\backslash u)]/u\}v$
\item $x[y\backslash \{[(yy)/z]u\}]=[\{x[y\backslash (yy)]\}/z]u$
\item $x[y\backslash \{[(yz)/y]u\}]=[\{x[y\backslash(yz)]\}/y]u$
\item $x[z\backslash \{[(yy)/y]u\}]=[\{x[z\backslash(yy)]\}/y]u$
\end{enumerate}
\end{myth}

\section{Main Results}
\subsection{Parastrophic invariance of groups}
\begin{myth}\label{0:12}
Let $G$ be a loop with identity element $e$ and $H$ a quasigroup
such that they are isotopic under the triple $\alpha =(A,B,C)$.
\begin{enumerate}
\item If $C=B$, then $G\overset{A}{\cong}H$ if and only if $eB\in N_\rho(H)$ where $N_\rho(H)$ is
the right nucleus of $H$.
\item If $C=A$, then $G\overset{B}{\cong}H$ if and only if $eA\in N_\lambda(H)$ where $N_\lambda(H)$ is
the left nucleus of $H$.
\end{enumerate}
\end{myth}
{\bf Proof}\\
Here, when $L_x$ and $R_x$ are respectively the left and right
translations of the loop $G$ then the left and right translations of
its quasigroup isotope $H$ are denoted by $L_x'$ and $R_x'$
respectively.

Let $(G,\cdot )$ and $(H,\circ )$ be any two distinct quasigroups.
If $A,B,C : G\rightarrow H$ are permutations, then the following
statements are equivalent :
\begin{itemize}
\item the triple $\alpha=(A,B,C)$ is an isotopism of $G$ upon $H$.
\item $R_{xB}'=A^{-1}R_xC~\forall~x\in G$.
\item $L_{yA}'=B^{-1}L_yC~\forall~y\in G$.
\end{itemize}
\begin{enumerate}
\item When $\alpha=(A,B,B)$, $R_{eB}'=A^{-1}B\Rightarrow
B=AR_{eB}'$. So,
\begin{displaymath}
\alpha=(A,AR_{eB}',AR_{eB}')=(A, A,
A)(I,R_{eB}',R_{eB}'),~\alpha~:~G\rightarrow H.
\end{displaymath}
If $(A, A, A,)~:~G\rightarrow H$ is an isotopism i.e $A$ is an
isomorphism, then $(I,R_{eB}',R_{eB}')~:~H\rightarrow H$ is an
autotopism if and only if $eB\in N_\rho(H)$.
\item When $\alpha=(A,B,A)$, $L_{eA}'=B^{-1}A\Rightarrow
A=BL_{eA}'$. So,
\begin{displaymath}
\alpha=(BL_{eA}',B,BL_{eA}')=(B,B,B)(L_{eA}',I,L_{eA}'),~\alpha~:~G\rightarrow
H.
\end{displaymath}
If $(B, B, B,)~:~G\rightarrow H$ is an isotopism i.e $B$ is an
isomorphism, then $(L_{eA}',I,L_{eA}')~:~H\rightarrow H$ is an
autotopism if and only if $eA\in N_\lambda(H)$.
\end{enumerate}

\begin{myth}\label{1:1}
A quasigroup $(S,\theta )$ is associative if and only if any of the
following equivalent statements is true.
\begin{enumerate}
\item $(S,\theta )$ is isotopic to $(S,(\theta^{-1})^*)$. Hence, the
other 4 parastrophes are also isotopic to $(S,(\theta^{-1})^*)$.
\item $(S,\theta^*)$ is isotopic to $(S,\theta^{-1})$. Hence, the
other 4 parastrophes are also isotopic to $(S,\theta^{-1})$.
\item $(S,\theta)$ is isotopic to $(S,({}^{-1}\theta )^*)$. Hence, the
other 4 parastrophes are also isotopic to $(S,({}^{-1}\theta )^*)$.
\item $(S,\theta^*)$ is isotopic to $(S,{}^{-1}\theta )$. Hence, the
other 4 parastrophes are also isotopic to $(S,{}^{-1}\theta )$.
\end{enumerate}
\end{myth}
{\bf Proof}\\
$(S,\theta )$ is associative if and only if $s_1\theta (s_2\theta
s_3)=(s_1\theta s_2)\theta s_3\Leftrightarrow
R_{s_2}R_{s_3}=R_{s_2\theta s_3}\Leftrightarrow L_{s_1\theta
s_2}=L_{s_2}L_{s_1}~\forall~s_1,s_2,s_3\in S$.

The proof of the equivalence of 1. and 2. is as follows.
$L_{s_1\theta s_2}=L_{s_2}L_{s_1}\Leftrightarrow {\cal L}_{s_1\theta
s_2}^{-1}={\cal L}_{s_2}^{-1}{\cal L}_{s_1}^{-1}\Leftrightarrow
{\cal L}_{s_1\theta s_2}={\cal L}_{s_1}{\cal L}_{s_2}\Leftrightarrow
(s_1\theta
s_2)\theta^{-1}s_3=s_2\theta^{-1}(s_1\theta^{-1}s_3)\Leftrightarrow
(s_1\theta s_2){\cal R}_{s_3}=s_2\theta^{-1}s_1{\cal
R}_{s_3}=s_1{\cal R}_{s_3}(\theta^{-1})^*s_2\Leftrightarrow
(s_1\theta s_2){\cal R}_{s_3}=s_1{\cal
R}_{s_3}(\theta^{-1})^*s_2\Leftrightarrow (s_2\theta^* s_1){\cal
R}_{s_3}=s_2\theta^{-1}s_1{\cal R}_{s_3}\Leftrightarrow $
\begin{equation}\label{eq:1}
({\cal R}_{s_3}, I, {\cal R}_{s_3})~:~(S,\theta )\to
(S,(\theta^{-1})^*)\Leftrightarrow
\end{equation}
\begin{equation}\label{eq:2}
(I, {\cal R}_{s_3},{\cal R}_{s_3})~:~(S,\theta^* )\to
(S,\theta^{-1})
\end{equation}
$\Leftrightarrow (S,\theta )$ is isotopic to
$(S,(\theta^{-1})^*)\Leftrightarrow (S,\theta^*)$ is isotopic to
$(S,\theta^{-1})$.

The proof of the equivalence of 3. and 4. is as follows.
$R_{s_2}R_{s_3}=R_{s_2\theta s_3}\Leftrightarrow
{\R}_{s_2}^{-1}{\R}_{s_3}^{-1}={\R}_{s_2\theta
s_3}^{-1}\Leftrightarrow {\R}_{s_3}{\R}_{s_2}={\R}_{s_2\theta
s_3}\Leftrightarrow (s_1{}^{-1}\theta s_3){}^{-1}\theta
s_2=s_1{}^{-1}\theta (s_2\theta s_3)\Leftrightarrow (s_2\theta
s_3){\tp}_{s_1}=s_3{\tp}_{s_1}{}^{-1}\theta s_2=s_2({}^{-1}\theta
)^*s_3{\tp}_{s_1}\Leftrightarrow (s_2\theta
s_3){\tp}_{s_1}=s_2({}^{-1}\theta )^*s_3{\tp}_{s_1}\Leftrightarrow
(s_3\theta^* s_2){\tp}_{s_1}=s_3{\tp}_{s_1}{}^{-1}\theta
s_2\Leftrightarrow$
\begin{equation}\label{eq:3}
(I, {\tp}_{s_1}, {\tp}_{s_1})~:~(S,\theta )\to
(S,({}^{-1}\theta)^*)\Leftrightarrow
\end{equation}
\begin{equation}\label{eq:4}
({\tp}_{s_1},I, {\tp}_{s_1})~:~(S,\theta^* )\to (S,{}^{-1}\theta
)
\end{equation}
$\Leftrightarrow (S,\theta)$ is isotopic to $(S,({}^{-1}\theta
)^*)\Leftrightarrow (S,\theta^*)$ is isotopic to $(S,{}^{-1}\theta
)$.

The last part of 1. to 4. follow by Theorem~\ref{0:1.3}.

\begin{myrem}
In the proof of Theorem~\ref{1:1}, it can be observed that the
isotopisms are triples of the forms $(A,I,A)$ and $(I,B,B)$. Weak
associative identities such as the Bol, Moufang and extra identities
have been found to be isotopic invariant in loops for any triple of
the form $(A,B,C)$ while in \cite{top}, the central identities have
been found to be isotopic invariant only under triples of the forms
$(A,B,A)$ and $(A,B,B)$. Since associativity obeys all the
Bol-Moufang identities, the observation in the theorem agrees with
the latter stated facts.
\end{myrem}

\begin{mycor}\label{1:2}
$(S,\theta )$ is an associative quasigroup if and only if any one of
particular four parastrophes of $(S,\theta )$ obeys any of the six
Khalil conditions or Evans law or a balanced identity or a
generalized group identity.
\end{mycor}
{\bf Proof}\\
Let $(S,\theta )$ be the quasigroup in consideration. By hypothesis,
$(S, \theta )$ is a group. Notice that $R_{s_2}R_{s_3}=R_{s_2\theta
s_3}\Leftrightarrow L_{s_2\theta s_3}^*=L_{s_3}^*L_{s_2}^*$. Hence,
$(S,\theta^*)$ is also a group. In Theorem~\ref{1:1}, two of the
parastrophes are isotopes of $(S,\theta )$ while the other two are
isotopes of $(S,\theta^*)$. Since the Khalil conditions, Evans law,
balanced identity and generalized group identity of
Theorem~\ref{0:1.2}, Theorem~\ref{0:1.4}, Theorem~\ref{0:1.5} and
Theorem~\ref{0:1.6} respectively, are necessary and sufficient
conditions for a quasigroup to be an isotope of a group, then they
must be necessarily and sufficiently true in the four quasigroup
parastrophes of $(S,\theta )$.

\begin{mycor}\label{1:21}
Let $(S,\theta )$ be an associative quasigroup.
\begin{enumerate}
\item $(S,\theta )\equiv (S,(\theta^{-1})^*)$ if and only if $(S,(\theta^{-1})^*)$ is associative.
\item $(S,\theta^*)\equiv (S,\theta^{-1})$ if and only if $(S,\theta^{-1})$ is associative.
\item $(S,\theta)\equiv (S,({}^{-1}\theta )^*)$ if and only if $(S,({}^{-1}\theta )^*)$ is associative.
\item $(S,\theta^*)\equiv (S,{}^{-1}\theta )$ if and only if $(S,{}^{-1}\theta )$ is associative.
\end{enumerate}
\end{mycor}
{\bf Proof}\\
Let $(S,\theta )$ be an associative quasigroup with identity element
$e$.
\begin{enumerate}
\item Using isotopism (\ref{eq:1}) of Theorem~\ref{1:1} in the second part of Theorem~\ref{0:12}, it will be observed that
$(S,\theta )\overset{I}{\cong}(S,(\theta^{-1})^*)$ if and only if
$s(\theta^{-1})^*e=s\in N_\lambda(S,(\theta^{-1})^*)~\forall~s\in S$
if and only if $(S,(\theta^{-1})^*)$ is associative.
\item Using isotopism (\ref{eq:2}) of Theorem~\ref{1:1} in the first part of Theorem~\ref{0:12}, it will be observed that
$(S,\theta^*)\overset{I}{\cong}(S,\theta^{-1})$ if and only if
$e\theta^{-1} s=s\in N_\rho(S,\theta^{-1})~\forall~s\in S$ if and
only if $(S,\theta^{-1})$ is associative.
\item Using isotopism (\ref{eq:3}) of Theorem~\ref{1:1} in the first part of Theorem~\ref{0:12}, it will be observed that
$(S,\theta)\overset{I}{\cong}(S,({}^{-1}\theta )^*)$ if and only if
$e({}^{-1}\theta )^*s=s\in N_\rho(S,({}^{-1}\theta )^*)~\forall~s\in
S$ if and only if $(S,({}^{-1}\theta )^*)$ is associative.
\item Using isotopism (\ref{eq:4}) of Theorem~\ref{1:1} in the second part of Theorem~\ref{0:12}, it will be observed that
$(S,\theta^*)\overset{I}{\cong}(S,{}^{-1}\theta )$ if and only if
$s({}^{-1}\theta )e=s\in N_\lambda(S,{}^{-1}\theta )~\forall~s\in S$
if and only if $(S,{}^{-1}\theta )$ is associative.
\end{enumerate}

\begin{mycor}\label{1:14}
Isotopy-isomorphy is a necessary and sufficient condition for any
two distinct quasigroups $(S,\pi_i)$ and $(S,\pi_j)$, $i,j\in
\{1,2,3,4,5,6\}$ to be parastrophic invariant under the associative
law.
\end{mycor}
{\bf Proof}\\
By Theorem~\ref{1:1} and Corollary~\ref{1:21}, the claim follows.

\subsection{Parastrophy-Holomorphy and Holomorphy-Parastrophy of Quasigroups}
\begin{myth}\label{2:5}
Let $(L,\theta)$ be a quasigroup with holomorph
$H(L,\theta)=(H,\circ)$. The following are true.
\begin{enumerate}
\item $(H,\circ^*)\equiv \big(H,(\circ_*)^*\big)\Leftrightarrow (L,\theta)\equiv (L,\theta^*)$.
\item $(H,\circ^{-1})\equiv \big(H,(\circ_{-1})^{-1}\big)\Leftrightarrow (L,\theta)\equiv (L,\theta^{-1})$.
\item $(H,{}^{-1}\circ)\equiv \big(H,{}^{-1}({}_{-1}\circ)\big)\Leftrightarrow (L,\theta)\equiv (L,{}^{-1}\theta)$.
\item $(H,(\circ^{-1})^*)\equiv \Big(H,(((\circ_{-1})_*)^{-1})^*\Big)\Leftrightarrow (L,\theta)\equiv (L,(\theta^{-1})^*)$.
\item $(H,({}^{-1}\circ )^*)\equiv \Big(H,({}^{-1}(({}_{-1}\circ)_*))^*\Big)\Leftrightarrow (L,\theta)\equiv (L,({}^{-1}\theta )^*)$.
\end{enumerate}
\end{myth}
{\bf Proof}\\
$(H,\circ)~:~(\alpha, x)\circ (\beta, y)=(\alpha\beta, x\beta\theta
y)$.
\begin{enumerate}
\item $(H,\circ^*)~:~(\beta, y)\circ (\alpha, x)=(\alpha\beta, x\beta\theta y)
=(\alpha\beta,y\theta^*x\beta)$. $(H,\circ_*)~:~(\alpha, x)\circ_*
(\beta, y)=(\alpha\beta, x\beta\theta^* y)=(\alpha\beta,y\theta
x\beta)\Rightarrow (\beta, y)(\circ_*)^*(\alpha, x)=(\alpha\beta,
x\beta\theta^* y)=(\alpha\beta,y\theta x\beta)$. Hence,
$(H,\circ^*)\equiv \big(H,(\circ_*)^*\big)\Leftrightarrow
(L,\theta)\equiv (L,\theta^*)$.
\item $(H,\circ^{-1})~:~(\alpha, x)\circ^{-1}(\alpha\beta,
x\beta\theta y)=(\beta, y)$. $(H,\circ_{-1})~:~(\alpha,
x)\circ_{-1}(\beta, y)=(\alpha\beta, x\beta\theta^{-1} y)\Rightarrow
(\alpha, x)(\circ_{-1})^{-1} (\alpha\beta, x\beta\theta^{-1}
y)=(\beta, y)$. Thence, $(H,\circ^{-1})\equiv
\big(H,(\circ_{-1})^{-1}\big)\Leftrightarrow (L,\theta)\equiv
(L,\theta^{-1})$.
\item $(H,{}^{-1}\circ)~:~(\alpha\beta, x\beta\theta y){}^{-1}\circ (\beta, y)=(\alpha,
x)$. $(H,{}_{-1}\circ)~:~(\alpha, x){}_{-1}\circ (\beta,
y)=(\alpha\beta, x\beta {}^{-1}\theta y)\Rightarrow (\alpha\beta,
x\beta {}^{-1}\theta y){}^{-1}({}^{-1}\circ )(\beta, y)=(\alpha,
x)$. Whence, $(H,{}^{-1}\circ)\equiv
\big(H,{}^{-1}({}_{-1}\circ)\big)\Leftrightarrow (L,\theta)\equiv
(L,{}^{-1}\theta)$.
\item $\big(H,(\circ^{-1})^*,\big)~:~(\alpha\beta,
x\beta\theta y)(\circ^{-1})^*(\alpha, x)=(\beta, y)$.
$(H,(\circ_{-1})_*)~:~(\alpha, x)(\circ_{-1})_*(\beta,
y)=(\alpha\beta, x\beta (\theta^{-1})^* y)\Rightarrow (\alpha\beta,
x\beta (\theta^{-1})^* y)(((\circ_{-1})_*)^{-1})^*(\alpha,
x)=(\beta, y)$. Then, $(H,(\circ^{-1})^*)\equiv
\big(H,(((\circ_{-1})_*)^{-1})^*\big)\Leftrightarrow
(L,\theta)\equiv (L,(\theta^{-1})^*)$.
\item $\big(H,({}^{-1}\circ)^*\big)~:~(\beta, y)({}^{-1}\circ )^*(\alpha\beta, x\beta\theta y)=(\alpha,
x)$. $(H,({}_{-1}\circ)_*)~:~(\alpha, x)({}_{-1}\circ )_* (\beta,
y)=(\alpha\beta, x\beta ({}^{-1}\theta )^* y)\Rightarrow (\beta,
y)((({}_{-1}\circ )_* )^{-1})^* (\alpha\beta, x\beta\theta
y)=(\alpha, x)$. So, $(H,({}^{-1}\circ )^*)\equiv
\big(H,({}^{-1}(({}_{-1}\circ)_*))^*\big)\Leftrightarrow
(L,\theta)\equiv (L,({}^{-1}\theta )^*)$.
\end{enumerate}

\begin{mycor}\label{2:6}
Let the quasigroup $(L,\theta)$ be a group with holomorph
$H(L,\theta)=(H,\circ)$. The following are true.
\begin{enumerate}
\item $(H,\circ^{-1})\equiv \Big(H,((\circ_{-1})^{-1})^*\Big)$ if and only if $(L,\theta^{-1})$ is associative.
\item $(H,{}^{-1}\circ)\equiv \Big(H,({}^{-1}({}_{-1}\circ))^*\Big)$ if and only if $(L,{}^{-1}\theta )$ is associative.
\item $(H,(\circ^{-1})^*)\equiv \Big(H,(((\circ_{-1})_*)^{-1})^*\Big)$ if and only if $(L,(\theta^{-1})^*)$ is associative.
\item $(H,({}^{-1}\circ )^*)\equiv \Big(H,({}^{-1}(({}_{-1}\circ)_*))^*\Big)$ if and only if $(L,({}^{-1}\theta )^*)$ is associative.
\end{enumerate}
\end{mycor}
{\bf Proof}\\
1. to 4. are proved by applying 4. and 5. of Theorem~\ref{2:5} to 1.
to 4. of Corollary~\ref{1:21}.

\begin{mycor}\label{2:7}
A necessary and sufficient condition for any two distinct
quasigroups $(S,\pi_i)$ and $(S,\pi_j)$, $i,j\in \{1,2,3,4,5,6\}$ to
be parastrophic invariant under the associative law is either
\begin{enumerate}
\item if the $\pi_i$-parastrophe of $H$ is equivalent to the
$\pi_i$-parastrophe of the holomorph of the $\pi_i$-parastrophe of
$S$ or
\item if the $\pi_i$-parastrophe of $H$ is equivalent to the $\pi_k$-parastrophe of the
$\pi_i$-parastrophe of the holomorph of the $\pi_i$-parastrophe of
$S$, for a particular $k\in \{1,2,3,4,5,6\}$.
\end{enumerate}
\end{mycor}
{\bf Proof}\\
The proofs of 1. and 2. can be deduced from Theorem~\ref{2:5} and
Corollary~\ref{2:6} respectively.

\paragraph{Further Study}
It will be interesting to investigate parastrophic invariance in
quasigroups relative to weak associative laws such as the
Bol-Moufang identities and establish necessary and sufficient
condition(s) for them to be parastrophic invariant in quasigroups.

\end{document}